\documentclass{amsart}        

\usepackage{amsmath,amsthm}   
\usepackage{amssymb}          
\usepackage{enumerate}        
\usepackage{euscript}         

\usepackage[ps,matrix,arrow,curve]{xy}       
\xymatrixcolsep{1.7pc} 
\xymatrixrowsep{1.7pc}
\newdir{ >}{{}*!/-5pt/\dir{>}}               

\newenvironment{roenumerate}{\begin{enumerate}[(i)]}{\end{enumerate}}

\newtheorem{thm}{Theorem}[section]
\newtheorem{defn}[thm]{Definition}
\newtheorem{prop}[thm]{Proposition}

\newtheorem{cor}[thm]{Corollary}
\newtheorem{lemma}[thm]{Lemma}

\theoremstyle{definition}  

\newtheorem{note}[thm]{Note}
\newtheorem{axiom}[thm]{Axiom}

\newcommand{\ulp}{\textup{(}}
\newcommand{\urp}{\textup{)}}

\newcommand{\dfn}{\textbf} 
\newcommand{\mdfn}[1]{\dfn{\mathversion{bold}#1}} 

\DeclareMathOperator{\Hom}      {Hom}
\DeclareMathOperator{\Ext}      {Ext}
\DeclareMathOperator{\PExt}     {PExt}
\DeclareMathOperator{\Ind}      {Ind}
\DeclareMathOperator{\ob}       {ob}
\DeclareMathOperator{\image}    {image}

\newcommand{\cat}               {\EuScript}  
\newcommand{\cA}                {{\cat A}}
\newcommand{\cB}                {{\cat B}}
\newcommand{\cC}                {{\cat C}}

\newcommand{\cF}                {{\cat F}}

\newcommand{\cH}                {{\cat H}}

\newcommand{\cP}                {{\cat P}}
\newcommand{\cS}                {{\cat S}}

\newcommand{\cSP}               {\cS/\cP}
\newcommand{\Ab}                {{\cat Ab}}

\newcommand{\yoneda}{\mathrm{h}}

\newcommand{\quo}               {/}
\newcommand{\Zmp}               {\Z \quo p}
\newcommand{\Zmn}               {\Z \quo n}
\newcommand{\Zmpi}              {\Z \quo p^{\infty}}
\newcommand{\Zmpj}              {\Z \quo p^j}
\newcommand{\Zmpk}              {\Z \quo p^k}
\newcommand{\QZ}                {\Q/\Z}

\newcommand{\pinf}              {p^{\infty}\!}

\newcommand{\Bh}                {\widehat{B}}
\newcommand{\Bt}                {\widetilde{B}}

\newcommand{\invlim}            {\varprojlim}
\newcommand{\dirlim}            {\varinjlim}

\newcommand{\xra}               {\xrightarrow}
\newcommand{\ra}                {\xrightarrow{}}
\newcommand{\la}                {\xleftarrow{}}

\newcommand{\Z}                 {\mathbb{Z}}
\newcommand{\Zp}                {\mathbb{Z}_p}
\newcommand{\Q}                 {\mathbb{Q}}
\newcommand{\C}                 {\mathbb{C}}
\newcommand{\I}                 {\mathbb{I}}

\newcommand{\Oplus}             {\displaystyle \bigoplus}
\newcommand{\Smash}             {\wedge}

\DeclareRobustCommand{\bigWedge}{\bigvee}

\newcommand{\ie}                {\emph{i.e.,} }
\newcommand{\iso}               {\cong}
\newcommand{\st}                {\;|\;}

\newcommand{\ab}                {\underline{a}}
\newcommand{\bb}                {\underline{b}}
\newcommand{\pinv}              {{\scriptstyle \frac{1}{p}}}

\renewcommand{\:}{\colon}

\newcommand{\period}    {{\makebox[0pt][l]{\hspace{2pt} .}}}
\newcommand{\comma}     {{\makebox[0pt][l]{\hspace{2pt} ,}}}

\begin{document}

\title[Phantom Maps]{ Phantom Maps and Homology Theories }
\author{J. Daniel Christensen}
\thanks{The first author was partially supported by an NSF grant and
  an NSERC scholarship.} 
\address{Department of Mathematics\\ M.I.T.\\ Cambridge, MA 02139} 
\email{jdchrist@math.mit.edu}

\author{Neil P. Strickland}
\address{Trinity College\\ Cambridge CB2 1TQ\\ England}
\thanks{The second author was partially supported by an NSF grant.}
\email{neil@pmms.cam.ac.uk}



\begin{abstract}
 We study phantom maps and homology theories in a stable homotopy
 category $\cS$ via a certain Abelian category $\cA$.  We express the
 group $\cP(X,Y)$ of phantom maps $X\ra Y$ as an $\Ext$ group in
 $\cA$, and give conditions on $X$ or $Y$ which guarantee that it
 vanishes.  We also determine $\cP(X,HB)$.  We show that any composite
 of two phantom maps is zero, and use this to reduce Margolis's
 axiomatisation conjecture to an extension problem.  We show that a
 certain functor $\cS\ra\cA$ is the universal example of a homology
 theory with values in an AB 5 category and compare this with some
 results of Freyd.
\end{abstract}

\maketitle

\tableofcontents

\section{Introduction}

In this paper we collect together a number of results about the
homotopy category of spectra.  A central theme is the problem of
reconstructing this category from the category of finite spectra or
(what is almost equivalent) from the category of generalised homology
theories.  A central result (to be explained in more detail below) is
that the category of spectra is a non-split linear extension of the
category of homology theories by a certain square-zero ideal, the
ideal of phantom maps.

Many of our results hold not only for the category of spectra but also
for other categories with similar formal properties.  In
Section~\ref{se:axioms}, we give a list of axioms which are sufficient
for most of the theory.  Let $\cS$ be a category satisfying these
axioms, and $\cF$ the full subcategory of finite objects.  In
Section~\ref{se:homology} we study the category $\cA$ of additive
functors from $\cF$ to the category $\Ab$ of Abelian groups, with
emphasis on the homology theories.  We also study the functor
$\yoneda\:\cS\ra\cA$ that sends a spectrum $X$ to the homology theory
$\yoneda_X$ it represents.

In Section~\ref{se:phantoms}, we consider phantom maps: a map 
$f\:X\ra Y$ is called \dfn{phantom} if
$\yoneda_f\:\yoneda_X\ra\yoneda_Y$ is zero, and the group of phantom
maps from $X$ to $Y$ is written $\cP(X,Y)$.  In Section~\ref{se:mac},
we show how our results about phantoms give new evidence for a
conjectured axiomatic characterisation of the classical stable
homotopy category, due to Margolis.  
In Section~\ref{se:pc}, we analyse the groups $\cP(X,HA)$, where $X$
is an arbitrary spectrum and $HA$ is an Eilenberg--Mac\,Lane
spectrum. 
Finally, in Section~\ref{se:uht}, we show
that the functor $\yoneda\:\cS\ra\cA$ is the universal example of a
homology theory on $\cS$ with values in an Abelian category satisfying
Grothendieck's axiom AB 5.  We also compare this with Freyd's
construction of a universal example without the AB 5 condition, and
make some related remarks about pro-spectra and ind-spectra.  

We next give a more detailed summary of our main results.  First we
show that there are several characterisations of phantom maps.

\begin{prop}
 Let $f\:X\ra Y$ be a map of spectra.  Then the following
 conditions are equivalent:
 \begin{roenumerate}
 \item $f$ is phantom, \ie $\yoneda_X(W)\ra\yoneda_Y(W)$ is zero for
  each finite $W$.
 \item $H(f)\:H(X)\ra H(Y)$ is zero for each homology theory $H$.
 \item The composite $W\ra X\ra Y$ is zero for each finite spectrum
  $W$ and each map $W \ra X$.
 \item The composite $X\ra Y\ra IW$ is zero for each finite
  spectrum $W$ and each map $Y\ra IW$.  \ulp{}Here $IW$ denotes the
  Brown--Comenetz dual of $W$; see Section~\ref{se:homology}.\urp{}
 \end{roenumerate}
\end{prop}

Another important result is the following.
\begin{thm} \label{th:composite}
 The composite of two phantom maps is zero \ulp{}and thus the phantom
 maps form a square-zero ideal\urp{}.
\end{thm}

This is a result that is folklore, but as far as we are aware the only
proof that works in this generality is the one presented here, which
was independently discovered by Neeman~\cite{ne:tba}.  Neeman also
proved some parts of Propositions~\ref{pr:1.4}, \ref{pr:1.5}
and~\ref{pr:1.6}.  Ohkawa~\cite{oh:vtaybkss} has a proof of
Theorem~\ref{th:composite} which works in the stable homotopy category
and uses CW-structures; it is not clear whether it goes through under
our axiomatic assumptions.  A simpler proof that works for a stricter
notion of phantom map appears in Gray's thesis~\cite{gr:oph} and is
published in~\cite{grmc:upm}.  The two notions coincide when the
source has finite skeleta.

It turns out that a number of interesting concepts can be described in
terms of the homological algebra of the Abelian category $\cA$.  As
usual, an object $F$ of $\cA$ is said to be \dfn{projective} if maps
from $F$ lift over epimorphisms, and \dfn{injective} if maps to $F$
extend over monomorphisms.  A spectrum $X$ is \mdfn{$\cA$-projective}
if $\yoneda_X$ is projective in $\cA$ and \mdfn{$\cA$-injective} if
$\yoneda_X$ is injective in $\cA$.

Here are two of our main results.

\begin{thm}
 There is a natural isomorphism
 $\cP(\Sigma^{-1} X,Y)\iso\Ext_\cA(\yoneda_X,\yoneda_Y)$.
\end{thm}

\begin{prop} \label{pr:1.4}
 Let $F \in \cA$.
 The following are equivalent:
 \begin{roenumerate}
 \item $F$ has finite projective dimension.
 \item $F$ has projective dimension at most one.
 \item $F$ is a homology theory.
 \item $F$ has injective dimension at most one.
 \item $F$ has finite injective dimension.
 \end{roenumerate}
\end{prop}

In view of the above, if an object $F$ of $\cA$ is projective or
injective, then it has the form $\yoneda_X$ for some spectrum $X$
(which is unique up to isomorphism).  The following result
describes those $X$ for which $\yoneda_X$ is projective or injective.

\begin{prop} \label{pr:1.5}
 Let $X$ be a spectrum.  Then the following are equivalent:
 \begin{roenumerate}
  \item $X$ is $\cA$-projective.
  \item $X$ is a retract of a wedge of finite spectra.
  \item $\cP(X,Y)=0$ for each spectrum $Y$.
 \end{roenumerate}
 Similarly, the following are equivalent:
 \begin{roenumerate}
  \item $X$ is $\cA$-injective.
  \item $X$ is a retract of a product of Brown--Comenetz duals of
   finite spectra.
  \item $\cP(Y,X)=0$ for each spectrum $Y$.
 \end{roenumerate}
\end{prop}

We also prove the following facts:
\begin{prop} \label{pr:1.6}
 \begin{enumerate}
 \item The category $\cA$ has enough injectives and projectives.
 \item Any spectrum $X$ sits in a cofibre sequence
  $P\ra Q\ra X\ra\Sigma P$, where $P$ and $Q$ are $\cA$-projective and 
  $X\ra\Sigma P$ is phantom.  The sequence
  $\yoneda_P\ra\yoneda_Q\ra\yoneda_X$ is a short exact sequence in
  $\cA$.  The map $X\ra\Sigma P$ is weakly initial among phantom maps
  out of $X$.
 \item Dually, any $X$ sits in a cofibre sequence
  $\Sigma^{-1}K\ra X\ra J\ra K$, where $J$ and $K$ are $\cA$-injective
  and $\Sigma^{-1}K\ra X$ is phantom.  The sequence
  $\yoneda_X\ra\yoneda_J\ra\yoneda_K$ is a short exact sequence in
  $\cA$.  The map $\Sigma^{-1}K\ra X$ is weakly terminal among phantom
  maps into $X$.
 \item $IX$ is $\cA$-injective for each $X$.  
 \item If $\pi_iY$ is finite for each $i$, then $Y$ is $\cA$-injective.
 \item If $\pi_iY$ is finitely generated for each $i$, then $\cP(X,Y)$
  is divisible for each $X$.
 \item The group $\cP(H\Zmp,Y)$ is always a vector space over $\Zmp$,
  and is nonzero \ulp{}and thus not divisible\urp{} for some $Y$.
 \item If $X$ is $\cA$-projective and $[X,W]=0$ for each finite $W$,
  then $X=0$.
 \end{enumerate}
\end{prop}

The above material appears in Sections~\ref{se:homology} and
\ref{se:phantoms}.  We warn the reader that while the results are for
the most part self-dual, the proofs are not.  

In Section~\ref{se:mac} we show how the stable homotopy category can
be viewed as a linear extension of the category of homology theories
by the bimodule of phantom maps.  Our point in making this rigorous is
that both the category of homology theories and the bimodule of
phantom maps are determined by the category of finite spectra, and so
we see that the category of spectra is determined up to extension by
the category of finite spectra.  Moreover, the goal of
Section~\ref{se:pc} is to prove that the extension is not split.  We
begin with the following result on phantom cohomology classes.  Here
$\PExt$ denotes the subgroup of $\Ext$ consisting of the pure or
phantom extensions, $HB$ denotes the Eilenberg--Mac\,Lane spectrum
with $\pi_0 HB = B$, and $H_*$ denotes integral homology.

\begin{thm}
  For any spectrum $X$ and Abelian group $B$ we have
  $\cP(X,HB) = \PExt(H_{-1} X, B)$.
\end{thm}

After submitting this paper, we discovered that this theorem is a
special case of some earlier results.  One such result is due to Huber
and Meier~\cite{hume:cticwc}.  They show that if $E_{*}(-)$ is a
homology theory of finite type, $B$ is an Abelian group, and
$F^{*}(-)$ is a cohomology theory fitting into a natural short exact
sequence
\[ 0 \ra \Ext(E_{n-1}(X),B) \ra F^{n}(X) \ra \Hom(E_{n}(X),B) \ra 0,
\]
then the subgroup of phantom cohomology classes in $F^{n}(X)$ is
isomorphic to $\PExt(E_{n-1}(X),B)$.  Taking $E = H$ and $F = HB$
gives our result.  Pezennec~\cite{pe:ptff} proves essentially the same
result, while Yosimura~\cite{yo:hcbpc} removes the finite type
hypothesis on the homology theory $E$ and concludes that the subgroup
of phantom cohomology classes is isomorphic to 
$\invlim^1 F^{n-1}(X_\alpha)$, where the $X_\alpha$ range over the
finite subspectra of $X$.  Ohkawa~\cite{oh:vtaybkss} also comes to
this conclusion, but without assuming the existence of $E$, $B$, or
the short exact sequence.  Another reference for this last result
is~\cite{ch:itcpgs}.

Using the theorem we are able to calculate all phantom maps between
Eilenberg--Mac\,Lane spectra.

\begin{cor}
 We have
 \[ \cP(\Sigma^k HA,HB) = \begin{cases}
     \PExt(A,B) & \text{ if } k = -1 \\
     0          & \text{ otherwise. }
    \end{cases}
 \]
\end{cor}

When we take $A = \Zmpi$ and $B = \bigoplus_k\Zmpk$ we can use the
above result and an explicit calculation to show that the phantom
sequence
\[ 0 \ra \cP(\Sigma^{-1}HA,HB) \ra \cS(\Sigma^{-1}HA,HB) \ra 
        \cA(\Sigma^{-1}HA,HB) \ra 0
\]
is not split.  This implies that the linear extension is also
not split.

In Section~\ref{se:uht} our main result is that $\yoneda\:\cS\ra\cA$
is the universal example of a homology theory with values in an AB 5
category.

\begin{prop}
 Let $\cC$ be an AB 5 category, and $K\:\cS\ra\cC$ a homology theory.
 Then there is an essentially unique strongly additive exact functor
 $K'\:\cA\ra\cC$ such that $K'\circ\yoneda\simeq K$.
\end{prop}

We also prove that the Ind completion of the category of finite
spectra is the category of homology theories.

We are indebted to Haynes Miller, Mike Hopkins, Mark Hovey, and the
rest of the MIT topology community not only for the many conversations
about the work presented here, but also for providing such a
stimulating environment.

\section{Axiomatic stable homotopy theory} \label{se:axioms}

Many of the properties of the stable homotopy category follow from a
collection of axioms which we state below.  These axioms are a slight
generalisation of those found in~\cite{ma:ssa}, and a specialisation
of those studied in~\cite{hopast:ash} (as one sees
using~\cite[Theorem~1.2.1]{hopast:ash}).  We shall say that an object
$X$ in an additive category $\cS$ is \dfn{small} if the functor
$\cS(X,-)$ preserves all coproducts that exist in $\cS$.

\begin{defn}\label{de:mbc}
 A \dfn{monogenic Brown category} is a category $\cS$ \ulp{}whose
 objects are called \dfn{spectra} and whose morphism sets are denoted
 $[-,-]$ or $\cS(-,-)$\urp{} satisfying the following axioms:
 \begin{enumerate}
 \item $\cS$ is triangulated \ulp{}and satisfies the octahedral
  axiom\urp{}.  Triangles are sometimes called cofibre sequences.
 \item $\cS$ has set-indexed coproducts.  The coproduct is usually
  written $\bigWedge$.
 \item $\cS$ is closed symmetric monoidal~\cite{ma:cwm}.  The
  multiplication is called the \dfn{smash product} and is denoted
  $\Smash$, the unit is denoted $S^0$, and the function spectra are
  denoted $F(X,Y)$.  The smash product and function spectrum functors
  are required to be compatible with coproducts and the triangulated
  structure, and all diagrams that one would expect to commute are
  required to.  See~\cite[Appendix~A]{hopast:ash} for more details.
 \item \label{it:small} $S^0$ is small.
 \item \label{it:gen} $S^0$ is a \mdfn{graded weak generator} for
  $\cS$: if $\pi_n X = 0$ for each $n \in \Z$ then $X = 0$, where
  $\pi_n X$ is defined to be $[S^n,X]$ and $S^n$ is $\Sigma^n S^0$.
 \item \label{it:rep} Homology theories and maps between them are
  representable --- see Section~\ref{se:homology} for an explanation
  of this axiom.
 \end{enumerate}
\end{defn}

\begin{note}
 If we replace axioms~\ref{it:small} and \ref{it:gen} with the weaker
 assumption that there exists a set of small graded weak generators,
 we get the notion of a \dfn{Brown category}.  Most, if not all, of
 what we discuss here goes through in this more general setting; we
 restrict ourselves to the monogenic setting only for simplicity.  In
 fact, one can get a long way without a symmetric monoidal structure.
\end{note}

The classical stable homotopy category, the derived category of a
countable commutative ring, the homotopy category of $G$-equivariant
spectra (for $G$ a compact Lie group) and suitable categories of
comodules over countable cocommutative Hopf algebras all form Brown
categories, the first two being monogenic.

An important subcategory of a monogenic Brown category $\cS$ is the
category $\cF$ of finite spectra which we define below.  Its
importance stems from the fact that a homology functor on $\cS$ is
determined by how it behaves on finite spectra.  Later, we will see
that even more of the structure of $\cS$ is captured by $\cF$.

We first make some auxiliary definitions.

\begin{defn}
 A \dfn{thick subcategory} $\cC$ of a triangulated category $\cS$ is a
 full subcategory which is closed under cofibres and retracts.  That
 is, if $X\ra Y\ra Z$ is a cofibre sequence with two of $X$, $Y$, and
 $Z$ in $\cC$, then so is the third; and if $X$ is in $\cC$ and $Y$ is
 a retract of $X$, then $Y$ is in $\cC$.  If $D$ is a class of spectra
 in $\cS$, then the \mdfn{thick subcategory generated by $D$} is the
 intersection of all thick subcategories containing $D$.
\end{defn}

The following definition was made and studied in~\cite{lemast:esht},
following work of Dold and Puppe.
\begin{defn}
 Write $DX=F(X,S^0)$.  A spectrum $X$ is \dfn{strongly dualizable}
 if the natural map $DX \Smash Y \ra F(X,Y)$ is an isomorphism for
 each $Y$.
\end{defn}

It is not hard to see that the following conditions on a spectrum $X$
are equivalent.  For a proof, see~\cite[Theorem 2.1.3]{hopast:ash}.
\begin{enumerate}
 \item $X$ lies in the thick subcategory generated by $S^0$.
 \item $X$ is small.
 \item $X$ is strongly dualizable.
\end{enumerate}

\begin{defn}
 We say that a spectrum $X$ is \dfn{finite} if it satisfies the above
 conditions, and we write $\cF$ for the category of finite spectra.
\end{defn}

One can show that $\cF$ has a small skeleton $\cF'$.  One can also
show that $\cF$ is closed under the functor $D$, and that there is a
natural map $X\ra D^2X$ that is an isomorphism when $X$ is finite, so
that $D$ gives an equivalence $\cF^{\textup{op}}\simeq\cF$.  We call
this equivalence \dfn{Spanier--Whitehead duality}.

In the case of the classical stable homotopy category, a spectrum is
finite if and only if it is a possibly desuspended suspension spectrum
of a finite CW-complex.

\section{Homology theories} 
\label{se:homology}

An additive functor from a triangulated category to an Abelian
category is \dfn{exact} if it sends cofibre sequences to exact sequences.  A
\dfn{homology theory} on a triangulated category $\cS$ is an exact
functor to an Abelian category which preserves the coproducts that
exist in $\cS$.  Unless we state otherwise, the target category will
always be taken to be the category $\Ab$ of Abelian groups.  It is
shown in~\cite[Section 4]{hopast:ash} that a homology theory defined
on $\cF$ has an essentially unique extension to a homology theory
defined on all of $\cS$, so the categories of homology theories on
$\cF$ and $\cS$ are equivalent.  More precisely, we have the following
result.
\begin{prop}\label{pr:lambda}
 For each spectrum $X$ there is a naturally defined small diagram
 $\Lambda(X)=\{X_\alpha\st\alpha\in A(X)\}$ of small spectra with
 compatible maps $X_\alpha\ra X$ such that for any homology theory $H$
 on $\cS$, the induced map $\dirlim_\alpha H(X_\alpha)\ra H(X)$ is an
 isomorphism.  Moreover, if $K$ is a homology theory defined on $\cF$
 and we define $\widehat{K}(X)=\dirlim_\alpha K(X_\alpha)$ then
 $\widehat{K}$ is the unique homology theory on $\cS$ extending $K$
 \ulp{}up to canonical isomorphism\urp{}.
\end{prop}
\begin{cor}\label{cor:lambda}
 If $W$ is finite then $[W,-]$ is a homology theory so
 \[ [W,X]=\dirlim_\alpha[W,X_\alpha]. \]
 In particular, we see that any map $W\ra X$ factors through some
 $X_\alpha$.
\end{cor}
In fact, if we take $\cF'$ to be a small skeleton of $\cF$, we can
define $A(X)$ to be the category of pairs $(U,u)$ where
$U\in\cF'$ and $u\:U\ra X$.  The diagram $\Lambda(X)$ is then just
the functor $A(X)\ra \cF$ sending $(U,u)$ to $U$.  

\begin{defn}
 The homology theory $\yoneda_X\:\cF\ra\Ab$ represented by a spectrum
 $X$ is the functor $\yoneda_X(W)=\pi_0(X\Smash W)$.  We shall write
 $\yoneda(X)$ instead of $\yoneda_X$ where this is typographically
 convenient.  We use the same symbol $\yoneda_X$ for the unique
 extension of this to a homology theory on all of $\cS$, which is
 again given by $\yoneda_X(W)=\pi_0(X\Smash W)=\yoneda_W(X)$.

 We also write $\cA$ for the Abelian category of additive functors
 from $\cF$ to $\Ab$.  This category has small Hom sets since $\cF$
 has a small skeleton.  Note that $\yoneda$ gives a functor
 $\cS\ra\cA$.   Note also that if $W$ is a finite spectrum then
 $\yoneda_W(Z)=[DW,Z]$; it follows easily that
 $[V,W]=\cA(\yoneda_V,\yoneda_W)$ when $V$ and $W$ are finite.
\end{defn}

We now give a more complete statement of Axiom~\ref{it:rep} of
Definition~\ref{de:mbc}.  This follows from the other axioms if
$\pi_*S^0$ is countable, but not otherwise.  (See~\cite{ne:tba} and
\cite[Section 4]{hopast:ash} for details.)

\begin{axiom}
 If $H$ is a homology theory on $\cF$ (taking values in $\Ab$), then
 there is a spectrum $Y$ in $\cS$ and a natural isomorphism
 $\yoneda_Y\ra H$.  Moreover, a natural transformation from
 $\yoneda_Y$ to $\yoneda_Z$ is always induced by a map from $Y$ to
 $Z$.  (This map need not be unique.  It turns out that a spectrum $Y$
 representing a given  homology theory is unique up to a non-unique
 isomorphism.)
\end{axiom}

\begin{note}
 A \dfn{cohomology theory} with values in an Abelian category $\cB$ is
 a homology theory with values in $\cB^{\text{op}}$.  It \emph{does}
 follow from the first five axioms that every cohomology theory on
 $\cS$ with values in $\Ab$ is of the form $[-,Y]$ for some $Y$.  By
 the Yoneda lemma, natural transformations are uniquely representable.
\end{note}

We record some basic facts about the functor $\yoneda$.
\begin{prop}
 The functor $\yoneda\:\cS\ra\cA$ preserves both products and
 coproducts, and it sends cofibre sequences to exact sequences.
\end{prop}
\begin{proof}
 It is easy to see that limits and colimits in a functor category such
 as $\cA$ are computed pointwise.  Thus, the first claim is that 
 \[ \yoneda(\prod_i X_i)(W) = \prod_i\yoneda(X_i)(W) \]
 for each small $W$.  This follows easily using
 $\yoneda(Y)(W)=[DW,Y]$.  The second claim is that
 \[ \yoneda(\bigWedge_i X_i)(W)=\bigoplus_i\yoneda(X_i)(W), \]
 which follows similarly using the smallness of $DW$.  The third claim
 is that for any cofibre sequence $X\ra Y\ra Z$, the resulting
 sequence 
 \[ \pi_0(X\Smash W) \ra \pi_0(Y\Smash W) \ra \pi_0(Z\Smash W) \]
 is exact, and this is clear.
\end{proof}

We can now start our study of homological algebra in the category
$\cA$. 

\begin{lemma} \label{le:proj}
 A finite spectrum $W$ is $\cA$-projective.  Hence, a retract of a
 wedge of finite spectra is $\cA$-projective.
\end{lemma}

\begin{proof}
 Let $W$ be a finite spectrum and suppose that
 $\alpha:\yoneda_W=[DW,-]\ra G$ is a natural transformation.  By the
 Yoneda Lemma it corresponds to an element of $GDW$.  If
 $\beta:F\ra G$ is an epimorphism then $FDW\ra GDW$ is as well, so
 $\alpha$ factors through $\beta$.  Thus $\yoneda_W$ is projective for
 $W$ finite.  But projectives are closed under coproducts and
 retracts, so if $X$ is a retract of a wedge of finite spectra, then
 $\yoneda_X$ is projective.
\end{proof}

We can use this to show that $\cA$ has enough projectives.

\begin{lemma}
 The category $\cA$ has enough projectives.
\end{lemma}

\begin{proof}
 Let $F\:\cF \ra \Ab$ be an additive functor, and choose a small
 skeleton $\cF'$ of $\cF$.  Then the natural map
 \[ \Oplus_{W \in \cF'}\; \Oplus_{\alpha \in FW} [W,-] \ra F \]
 is clearly an epimorphism.  Using the fact that $[W,-]=\yoneda_{DW}$,
 we see that the source is projective.
\end{proof}

Note that the source of the above epimorphism is just $\yoneda_X$,
where 
\[ X=\bigWedge_{W \in \cF'} \;\bigWedge_{\alpha \in FW} W. \]
Moreover, $\yoneda_X$ is not just projective, but free in the
following sense.  Let $\cC$ be the category of $\ob(\cF')$-indexed
families of sets, and consider the evident forgetful functor
$\cA\ra\cC$.  This has a left adjoint, whose image consists of the
functors $\yoneda_X$, where $X$ is a wedge of finite spectra; it is
natural to regard these as the free objects of $\cA$.  As usual, an
object is projective if and only if it is a retract of a free object;
it follows that projective objects are homology theories.

\begin{lemma} \label{le:yon}
 A map $f\:X\ra Y$ is an isomorphism if and only if
 $\yoneda_f\:\yoneda_X\ra\yoneda_Y$ is an isomorphism.  The same holds
 with ``isomorphism'' replaced by ``split monomorphism'' or ``split
 epimorphism''.
\end{lemma}

\begin{proof}
 Suppose that $\yoneda_f\:\yoneda_X\ra\yoneda_Y$ is an isomorphism.
 As $\pi_k(X)=h_X(S^{-k})$, we see that $\pi_*(f)$ is an isomorphism,
 so $f$ is an isomorphism.

 If $\yoneda_f\:\yoneda_X\ra\yoneda_Y$ is a split monomorphism, choose
 a splitting, which by Brown Representability is of the form
 $\yoneda_g$.  The composite $\yoneda_g\circ\yoneda_f$ is the identity,
 so $gf$ is an isomorphism.  By composing $g$ with the inverse of this
 isomorphism we get a splitting of $f$.

 The case when $\yoneda_f$ is a split epimorphism is dual.
\end{proof}

\begin{prop} \label{pr:proj-char}
 A spectrum $X$ is $\cA$-projective if and only if it is a retract of
 a wedge of finite spectra.
\end{prop}

\begin{proof}
 $\Leftarrow$: This is Lemma~\ref{le:proj}.

 $\Rightarrow$: If $\yoneda_X$ is projective, it is a retract of
 $\yoneda_Y$ with $Y$ a wedge of finite spectra.  By Brown
 Representability and the previous lemma, $X$ is a retract of $Y$.
\end{proof}

\subsection{The dual picture}

We first recall the basic facts about duality for Abelian groups.
\begin{defn}
 For any Abelian group $A$, we write $\I(A)=\Hom(A,\QZ)$.  It is
 well-known that this is a contravariant exact functor which converts
 sums to products, and that the natural map $A\xra{}\I^2(A)$ is a
 monomorphism.  Moreover, if $A$ is finitely generated then $\I^2(A)$
 is the profinite completion of $A$; in particular, if $A$ is finite
 then $\I^2(A)=A$.
\end{defn}

Given a spectrum $X$ consider the contravariant functor from $\cS$ to
$\Ab$ sending $Y$ to $\I(\pi_0(X \Smash Y))$; this is clearly a
cohomology theory.  There is thus a representing object $IX$ such that
$\I(\pi_0(X \Smash Y))\simeq [Y,IX]$; we call this the
\dfn{Brown--Comenetz dual} of $X$~\cite{brco:pdg}.

\begin{prop} \label{pr:IX-inj}
 For each spectrum $X$, $IX$ is $\cA$-injective.
\end{prop}

\begin{proof}
 Fix a spectrum $X$.  As in Corollary~\ref{cor:lambda}, we have a
 diagram $\{X_\alpha\}$ of finite spectra such that
 $[W,X]=\dirlim_\alpha[W,X_\alpha]$ for all finite $W$.  We
 temporarily write $\cA'$ for the category of contravariant additive
 functors from $\cF$ to $\Ab$.  If $F$ is in $\cA$ we have
 \begin{align*}
  \cA(F, \yoneda_{IX}) 
  &= \cA(F, \I[-,X])                    \\
  &= \cA'([-,X], \I F)                  \\
  &= \cA'(\dirlim [-,X_\alpha],\I F)    \\
  &= \invlim \cA'([-,X_\alpha],\I F)    \\
  &= \invlim \I FX_\alpha               \\
  &= \I(\dirlim FX_\alpha).
 \end{align*}
 Suppose now that $F\ra G$ is a monomorphism in $\cA$.  We must show
 that the map $\cA(F,\yoneda_{IX})\la\cA(G,\yoneda_{IX})$ is a
 surjection.  Each map $FX_\alpha\ra GX_\alpha$ is monic, and a
 filtered colimit of monomorphisms is monic, so the map
 $\I(\dirlim FX_\alpha)\la \I(\dirlim GX_\alpha)$ is surjective,
 since $\QZ$ is injective.  Thus
 $\cA(F,\yoneda_{IX})\la\cA(G,\yoneda_{IX})$ is surjective. 
\end{proof}

\begin{cor}
 If $Y$ has finite homotopy groups, then $Y\iso I^2Y$ and so $Y$ is
 $\cA$-injective.  Moreover, for any family $\{X_i\}$ of spectra, the
 product $\prod_iI(X_i)=I(\bigWedge_{\!i} X_i)$ is $\cA$-injective, as
 is any retract of such a product. \qed
\end{cor}

\begin{prop} \label{pr:enough-inj}
 $\cA$ has enough injectives.
\end{prop}

\begin{proof}
 For finite $W$ a natural transformation from $G\in\cA$ to
 $\yoneda_{IW}$ corresponds to an element of $\I G(W)$.  Let
 $\cF'$ be a small skeleton of $\cF$, so there is a natural map
 \[ G\ra
    \prod_{W\in\cF'}\;\prod_{\alpha\in\I G(W)}\yoneda_{IW}.
 \]
 Since $\QZ$ is an injective cogenerator in the category of Abelian
 groups, one can show that this map is a monomorphism.
\end{proof}

In fact, the target of the monomorphism is the homology theory
represented by a product of Brown--Comenetz duals of finite spectra.
In $\cA$, being injective is equivalent to being a retract of such a
functor.  In particular, injectives are homology theories.

\begin{prop}
 A spectrum $X$ is $\cA$-injective if and only if it is a retract of a
 product of Brown--Comenetz duals of finite spectra.
\end{prop}

\begin{proof}
 $\Leftarrow$:  This follows from Proposition~\ref{pr:IX-inj}.

 $\Rightarrow$: If $\yoneda_X$ is injective, it is a retract of
 $\yoneda(\prod IW_\alpha)$ with each $W_\alpha$ finite.  As in the
 proof of Proposition~\ref{pr:proj-char}, this implies that $X$ is a
 retract of $\prod IW_\alpha$.
\end{proof}

\section{Phantom maps} \label{se:phantoms}


There is a class of maps that we cannot see, at least not easily.

\begin{prop} \label{pr:ph1}
 The following conditions on a map $f\:X\ra Y$ are equivalent:
 \begin{roenumerate}
 \item The natural transformation $\yoneda_f\:\yoneda_X\ra\yoneda_Y$
  is zero.
 \item For each homology theory $H$, we have $H(f)=0$. \label{it:ph}
 \item The composite $W\ra X\ra Y$ is zero for each finite spectrum
  $W$ and each map $W \ra X$. 
 \end{roenumerate}
 \ulp{}A fourth equivalent condition appears in
 Proposition~\ref{pr:ph2}.\urp{}
\end{prop}

\begin{proof}
 (iii)$\Rightarrow$(ii): 
 Let $\Lambda(X)=\{X_\alpha\}$ be as in Proposition~\ref{pr:lambda},
 so that $H(X)=\dirlim_\alpha H(X_\alpha)$.  The composite
 $X_\alpha\ra X\ra Y$ is zero by (iii), so $H(X_\alpha)\ra H(Y)$ is
 zero.  It follows that $H(f)\:H(X)\ra H(Y)$ is zero.

 (ii)$\Rightarrow$(i):
 Suppose that $H(f)=0$ for each homology theory $H$.  Then for each
 finite spectrum $W$, the map
 $\yoneda_W(f)\:\pi_0(X\Smash W)\ra\pi_0(Y\Smash W)$ is zero.  In other
 words, the natural map $h_f$ is zero at $W$.

 (i)$\Rightarrow$(iii): Suppose that (i) holds and that $W$ is finite.
 Then $DW$ is also finite, and $f$ induces the zero map
 $[W,X]=\pi_0(DW\Smash X)\ra\pi_0(DW\Smash Y)=[W,Y]$.
\end{proof}

\begin{defn}
 A map $X \ra Y$ satisfying the equivalent conditions of the
 proposition is called \dfn{phantom} or \mdfn{$\cA$-null}.  The
 collection of phantom maps from $X$ to $Y$ is denoted $\cP(X,Y)$ and
 is a subgroup of $[X,Y]$.  Similarly, we say that a map $X\ra Y$ is
 \mdfn{$\cA$-monic} or \mdfn{$\cA$-epic} if the natural transformation
 $\yoneda_X\ra\yoneda_Y$ is monic or epic, respectively.
\end{defn}

If $\{X_\alpha\}$ is an indexed collection of spectra, then the map
$\bigWedge X_\alpha\ra\prod X_\alpha$ is $\cA$-monic, and hence its
fibre is phantom.  As an example of this in the classical stable
homotopy category, let $C$ be the cokernel of the map from the sum of
countably many copies of $\Z$ to the product.  The fibre of the map
$H(\bigoplus\Z)\ra H(\prod\Z)$ between Eilenberg--Mac\,Lane spectra is
a phantom map $\Sigma^{-1}HC\ra H(\bigoplus\Z)$.  It is non-zero
because the short exact sequence $0\ra\bigoplus\Z\ra\prod\Z\ra C\ra 0$
is not split.  To see that this sequence is not split, notice that the
coset of the quotient containing $(1,2,4,8,16,\ldots)$ is non-zero and
is divisible by $2^k$ for each $k$, since initial terms may be dropped
without changing the coset.  But $\prod\Z$ contains no such elements,
so $C$ could not be a summand.  We learned this argument from Dan
Dugger, who credits it to~\cite{gr:nth}.

As further evidence of the ubiquity of phantom maps, it can be shown
that in the classical stable homotopy category there are uncountably
many phantom maps from $\C P^\infty$ to $S^3$.  Gray~\cite{gr:ssn} has
a proof for spaces which simplifies when read stably.

\begin{note}\label{no:sp}
 It is not hard to see that phantom maps form an \dfn{ideal} in $\cS$:
 if $f$, $g$ and $h$ are composable and $g$ is phantom, then $fg$ and
 $gh$ are phantom; and if $f$ and $g$ are parallel phantom maps, then
 $f+g$ is phantom.  This means that there is a well-defined additive
 category $\cSP$ having the same objects as $\cS$ and with
 $\cSP(X,Y):=\cS(X,Y)/\cP(X,Y)$.  We have a natural isomorphism
 $\cA(\yoneda_X,\yoneda_Y)\iso\cSP(X,Y)$, so $\yoneda$ gives an
 equivalence between $\cSP$ and the category $\cH$ of homology
 theories.
\end{note}

\begin{lemma} \label{le:univ-phantom}
 For any spectrum $X$ there is a weakly initial phantom map
 \[ \delta\:X\ra\widetilde{X} \]
 from $X$.  By `weakly initial' we mean that any other phantom map
 from $X$ factors through $\delta$, but we don't insist upon
 uniqueness.
\end{lemma}

\begin{proof}
 Let $\Lambda(X)=\{X_\alpha\}$ be as in Proposition~\ref{pr:lambda}.
 For each $\alpha$ we have a given map $X_\alpha\ra X$, so we get a
 map $\bigWedge_\alpha X_\alpha\ra X$.  Let
 $\delta\:X\ra\widetilde{X}$ be the cofibre of this map.
 Corollary~\ref{cor:lambda} tells us that every map from a finite
 spectrum $W$ to $X$ factors through $\bigWedge_\alpha X_\alpha$, so
 the composite $W\ra X\xra{\delta}\widetilde{X}$ is zero.  It
 follows that $\delta$ is phantom.  Moreover, any phantom map from $X$
 is zero when restricted to $\bigWedge_\alpha X_\alpha$ and so factors
 through $\delta$.
\end{proof}

\begin{cor}
 Let $X$ be a spectrum.  Then the following are equivalent:
 \begin{roenumerate}
  \item $X$ is $\cA$-projective.
  \item $X$ is a retract of a wedge of finite spectra.
  \item $\cP(X,Y)=0$ for each spectrum $Y$.
 \end{roenumerate}
\end{cor}

\begin{proof}
 If there are no phantom maps from $X$, then the weakly initial phantom map
 $X\ra\tilde{X}$ is zero, and so $X$ is a retract of the wedge of
 finite spectra $\bigWedge X_\alpha$.  Conversely, if $X$ is a retract
 of a wedge of finite spectra, then it is clear that there are no
 phantoms from $X$.

 Proposition~\ref{pr:proj-char} tells us that being a retract of a
 wedge of finite spectra is equivalent to being $\cA$-projective.
\end{proof}

\begin{prop} \label{pr:cofibre}
 Any spectrum $X$ sits in a cofibre sequence $P\ra Q\ra X\ra\Sigma P$,
 where $P$ and $Q$ are $\cA$-projective and $X\ra\Sigma P$ is phantom.
 The sequence $\yoneda_P\ra\yoneda_Q\ra\yoneda_X$ is a short exact
 sequence in $\cA$.  The map $X\ra\Sigma P$ is weakly initial among
 phantom maps out of $X$.
\end{prop}

\begin{proof}
 Consider the diagram
 \[ \xy \xymatrix{
      \bigWedge_{\!\alpha \ra \beta} X_\alpha \ar[r]^-{1-s} 
          \ar@/l1em/[d]
    & \bigWedge_{\!\alpha} X_\alpha \ar@{=}[d] \ar[r]
    & Y \ar@/l1em/[d] \\
      P \ar[r] \ar@/r1em/[u] 
    & \bigWedge_{\!\alpha} X_\alpha \ar[r]
    & X \period \ar@/r1em/[u]
  }\endxy 
 \]
 The map called $1$ includes the $\alpha\ra\beta$ summand into the
 $\alpha$ summand via the identity map, while the map $s$ (for
 `shift') sends the $\alpha\ra\beta$ summand to the $\beta$ summand
 via the map $X_\alpha \ra X_\beta$.  The map
 $\bigWedge_\alpha X_\alpha \ra X$ is the map considered in
 Lemma~\ref{le:univ-phantom}.  The spectra $Y$ and $P$ are defined to
 make the rows cofibre sequences, so $\widetilde{X}$ (from the lemma) is
 $\Sigma P$.  
 The composite 
 $\bigWedge_{\alpha \ra \beta} X_\alpha \ra \bigWedge_\alpha X_\alpha \ra X$ 
 is null, so there is a map of cofibre sequences in the downward direction.
 Now consider the following natural transformation from $\yoneda_X$ to
 $\yoneda_Y$.
 Let $W$ be a finite spectrum.
 An element of $\yoneda_X(W)$ is a map $DW \ra X$.  
 $DW$ is finite, so this map is $X_\gamma \ra X$ for some $\gamma$.
 We have a map $\bigWedge_\alpha X_\alpha \ra Y$, so in particular we have a map
 $X_\gamma \ra Y$.  That is, we have a map $DW \ra Y$, or an element
 of $\yoneda_Y(W)$.
 This defines a natural transformation $\yoneda_X \ra \yoneda_Y$,
 and by Brown Representability this natural transformation is
 induced by a map $X \ra Y$.
 By definition, the square commutes up to phantoms, but since
 $\bigWedge_\alpha X_\alpha$ is $\cA$-projective, the square commutes.
 One thus obtains a fill-in map $P \ra \bigWedge_{\alpha \ra \beta} X_\alpha$.
 Also, one can check that the composite $X \ra Y \ra X$ is an isomorphism,
 and it follows that the composite 
 $P \ra \bigWedge_{\alpha \ra \beta} X_\alpha \ra P$ is an isomorphism as well.
 Thus, $P$ is a retract of a wedge of finite spectra,
 and we have demonstrated that $X$ is the cofibre of a map between
 $\cA$-projective spectra.  We saw in Lemma~\ref{le:univ-phantom} that
 the map $X\ra\Sigma P$ is weakly initial.
\end{proof}

We now get an easy proof of a result that is folklore.  
The method of
proof presented in this section was independently discovered by
Neeman~\cite{ne:tba}.  
A proof for the special case of the classical stable homotopy category
was given by Ohkawa~\cite{oh:vtaybkss}.
A proof assuming that the source has finite
skeleta appears in \cite{gr:oph} and \cite{grmc:upm}.  
(See the introduction for more detailed comments.)

\begin{cor}
 The composite of two phantom maps is zero.
\end{cor}

\begin{proof}
 Suppose that $X \xra{f} Y$ and $Y \xra{g} Z$ are phantom.  Factor
 $f$ through $\delta$:
 \[ \xy \xymatrix{
   X \ar[r]^-{\delta} \ar[d]^-{f}
 & \Sigma P \ar@{{}-->}[dl]^-{f'} \\
   Y \ar[d]^-{g} \\
   Z \period } \endxy 
 \]
 The $\cA$-projectivity of $\Sigma P$ implies that $g f' = 0$ and so
 $g f = 0$.
\end{proof}

We can now characterise homology theories in terms of the homological
algebra of the category $\cA$.

\begin{prop} \label{pr:proj-1}
 A functor in $\cA$ is a homology theory if and only if it has finite
 projective dimension if and only if it has projective dimension at
 most one.
\end{prop}

\begin{proof}
 First, consider a short exact sequence $F\ra G\ra H$ in $\cA$, in
 which two of $F$, $G$ and $H$ are homology theories.  We claim that
 the third is also.  Indeed, consider a cofibre sequence
 $X\ra Y\ra Z$.  By applying $F$, we get a chain complex
 \[ \cdots \ra F(\Sigma^{-1}Z) \ra FX \ra FY \ra FZ \ra
     F(\Sigma X) \ra \cdots.
 \]
 By doing the same with $G$ and $H$, we obtain a short exact sequence
 of chain complexes.  By assumption, two of the three chain complexes
 are exact; it follows easily that the third is also, as required.

 We have seen that projective functors are homology theories.  It
 follows easily from the above that functors of finite projective
 dimension are homology theories (by induction on dimension).

 Consider a homology theory $H$.  There exists a spectrum $X$ such
 that $H=\yoneda_X$, and a cofibre sequence $P\ra Q\ra X\ra\Sigma P$
 as in Proposition~\ref{pr:cofibre}.  This gives a projective
 resolution $0\ra\yoneda_P\ra\yoneda_Q\ra\yoneda_X=H\ra 0$, so $H$ has
 projective dimension at most one.
\end{proof}

We can now describe the phantom maps in terms of $\cA$.

\begin{thm} \label{th:pff}
 The group $\cP(\Sigma^{-1} X, Y)$ of phantom maps is naturally
 isomorphic to $\Ext_\cA(\yoneda_X, \yoneda_Y)$.
\end{thm}

\begin{proof}
 Consider the usual projective resolution
 $0\ra\yoneda_P\ra\yoneda_Q\ra\yoneda_X\ra 0$ of $\yoneda_X$ in $\cA$.
 The first cohomology group of the left column of
 \[ \xy \xymatrix{
   0 &
   0 \\
   \cA(\yoneda_P, \yoneda_Y) \ar[u] \ar@{=}[r] &
   [P,Y] \ar[u] \\
   \cA(\yoneda_Q, \yoneda_Y) \ar[u] \ar@{=}[r] &
   [Q,Y] \ar[u] \\
   0 \ar[u] &
   0 \ar[u] }
 \endxy 
 \]
 is the $\Ext$ group in question, and the left column can be
 identified with the right column since $P$ and $Q$ are
 $\cA$-projective.  But the first cohomology of the right column is
 $\cP(\Sigma^{-1} X, Y)$ because every phantom $\Sigma^{-1} X \ra Y$
 extends to $P$, and the difference between two such extensions
 factors through $Q$.

 It is easy to see that the isomorphism is natural in $X$ and $Y$.
\end{proof}

\begin{note}
 The above proposition can also be proved using the definition of
 $\Ext$ in terms of equivalence classes of short exact sequences.  The
 isomorphism sends a phantom map $f\:\Sigma^{-1}X\ra Y$ to the short
 exact sequence
 \[ 0\ra\yoneda(Y)\ra\yoneda(\text{cofibre} f)\ra\yoneda(X)\ra 0.\]
\end{note}

\subsection{The dual picture}

Now we prove the dual results, making use of what came above.

\begin{prop} \label{pr:IY}
 For any spectra $X$ and $Y$, we have $\cP(X,IY) = 0$.
\end{prop}

\begin{proof}
 By Theorem~\ref{th:pff}, $\cP(\Sigma^{-1} X,IY) =
 \Ext_\cA(\yoneda_X,\yoneda_{IY})$.  But $\yoneda_{IY}$ is injective,
 so this is zero.

 One can prove this directly from the definition of $IY$ as well.
\end{proof}

With this we can now prove our fourth characterisation of phantom maps.

\begin{prop} \label{pr:ph2}
 A map $X \ra Y$ is phantom if and only if the composite
 $X\ra Y\ra IW$ is null for each finite $W$ and each map $Y \ra IW$.
\end{prop}

\begin{proof}
 By the previous proposition, every phantom map is null when composed
 with a map $Y \ra IW$.  

 Conversely, suppose that $X\ra Y$ is such that $(X\ra Y\ra IW)=0$ for
 all $Y\ra IW$.  Consider the spectrum
 \[ Z = \prod_{W \in \cF'} \; \prod_{Y \ra IW} IW. \]
 The evident map $Y\ra Z$ is $\cA$-monic, as in
 Proposition~\ref{pr:enough-inj}.  Since
 \[ \yoneda_X \ra \yoneda_Y \ra \yoneda_Z \]
 is null by assumption, the map $\yoneda_X\ra\yoneda_Y$ must also be
 null, so $X\ra Y$ is phantom.
\end{proof}

\begin{lemma} \label{le:BC-monic}
 There is a natural map $X\ra I^2 X$, which is $\cA$-monic for all
 $X$.
\end{lemma}

\begin{proof}
 Consider $[X,I^2 X]$.  By using the definition of $I$ twice, we find
 $[X,I^2 X] = \I(\pi_0(X \Smash IX)) = [IX,IX]$, and so there
 is a natural map $X \ra I^2 X$ corresponding to the identity map
 in $[IX,IX]$.

 We need to show that for $W$ finite, the map $[W,X]\ra [W,I^2 X]$ is
 monic.  We can calculate the latter group and we find that it is
 $\I^2[W,X]$.  The map $[W,X]\ra [W,I^2 X]$ is the natural inclusion
 of $[W,X]$ into its double dual; since $\QZ$ is an injective
 cogenerator, this is monic.
\end{proof}

\begin{prop} \label{pr:cofibre-inj}
 Any spectrum $X$ sits in a cofibre sequence
 $\Sigma^{-1}K\ra X\ra J\ra K$, where $J$ and $K$ are $\cA$-injective
 and $\Sigma^{-1}K\ra X$ is phantom.  The sequence
 $\yoneda_X\ra\yoneda_J\ra\yoneda_K$ is a short exact sequence in
 $\cA$.  The map $\Sigma^{-1}K\ra X$ is weakly terminal among phantom
 maps into $X$.
\end{prop}

\begin{proof}
 Let $J=I^2X$ and let $K$ be the cofibre of the natural map
 $X\ra I^2X$.  Similarly, let $L=I^2K$ and form the cofibre sequence
 $K\ra L\ra M$.  By Lemma~\ref{le:BC-monic} the maps
 $\Sigma^{-1}K\ra X$ and $\Sigma^{-1}M\ra K$ are phantom and so the
 cofibre sequences $X\ra J\ra K$ and $K\ra L\ra M$ become short exact
 in $\cA$.  Thus
 $\Ext_\cA^1(\yoneda_M,\yoneda_K) = \Ext_\cA^2(\yoneda_M,\yoneda_X)$,
 which vanishes as $\yoneda_M$ has projective dimension at most one
 (Proposition~\ref{pr:proj-1}).  Therefore the extension
 $\yoneda_K\ra\yoneda_L\ra\yoneda_M$ splits in $\cA$ and hence
 $\yoneda_K$ is injective.

 We showed in Proposition~\ref{pr:IY} that $\cP(-,J)=0$, and it
 follows easily that $\Sigma^{-1}K\ra X$ is weakly terminal.
\end{proof}

The following corollary is an easy consequence of the above
constructions. 
\begin{cor}
 The following are equivalent:
 \begin{roenumerate}
  \item $X$ is $\cA$-injective.
  \item $X$ is a retract of a product of Brown--Comenetz duals of
   finite spectra.
  \item $\cP(Y,X)=0$ for each spectrum $Y$.\qed
 \end{roenumerate}
\end{cor}

For example, this means that the completed Johnson-Wilson spectrum
$\widehat{E(n)}$ is $\cA$-injective.  Indeed, if $W$ is finite then
$\widehat{E(n)}^*W$ is compact Hausdorff in the $I_n$-adic topology.
The inverse limit functor is exact for inverse systems of compact
Hausdorff topological groups, and one can deduce from this that there
are no phantom maps to $\widehat{E(n)}$.

Summarising our homological results gives:

\begin{thm}\label{th:hpi}
 Let $F \in \cA$.  Then the following are equivalent:
 \begin{roenumerate}
 \item $F$ has finite projective dimension.
 \item $F$ has projective dimension at most one.
 \item $F$ is a homology theory.
 \item $F$ is in the image of $\yoneda$.
 \item $F$ has finite injective dimension.
 \item $F$ has injective dimension at most one.  \qed
 \end{roenumerate}
\end{thm}

\subsection{Divisibility}

To start with, we recall a result that is well-known to the experts.

\begin{prop} \label{pr:divisible}
 If $\pi_i Y$ is a finitely generated Abelian group for each $i$, then
 $\cP(X,Y)$ is divisible for each $X$.
\end{prop}

\begin{proof}
 Let $Z$ be the cofibre of the natural map $Y\ra I^2Y$.  We have seen
 that the resulting map $\Sigma^{-1}Z\ra Y$ is a weakly terminal
 phantom map, so that $\cP(X,Y)$ is a quotient of $[X,\Sigma^{-1}Z]$.  It
 will thus be enough to show that $[X,\Sigma^{-1}Z]$ is a rational
 vector space.

 The induced map $\pi_k(Y)\ra\pi_k(I^2Y)$ is just the inclusion of
 $\pi_k(Y)$ into its double dual with respect to $\QZ$, which is the
 same as its profinite completion (as $\pi_k(Y)$ is finitely
 generated).  It follows that $\pi_k(Z)$ is a finite direct sum of
 copies of $\widehat{\Z}/\Z$, which is well-known to be a rational
 vector space.  It follows that any nonzero integer $n$ induces an
 isomorphism $\pi_*(Z)\ra\pi_*(Z)$, and thus an isomorphism $Z\ra Z$.
 It follows that $[X,\Sigma^{-1}Z]$ is a rational vector space, as
 required.  
\end{proof}

It is not the case that $\cP(X,Y)$ is always divisible, however.
Indeed, we have the following result.
\begin{prop}\label{pr:nondiv}
 Let $\cS$ be the classical stable homotopy category, and $H\Zmp$ the
 mod $p$ Eilenberg--Mac\,Lane spectrum in $\cS$.  Then $\cP(H\Zmp,Y)$
 is a vector space over $\Zmp$, and there exist spectra $Y$ for which
 it is nonzero \ulp{}and thus not divisible\urp{}.
\end{prop}
\begin{proof}
 As $p$ times the identity map of $H\Zmp$ is zero, we see that
 $[H\Zmp,Y]$ is a vector space over $\Zmp$, so the same is true of
 $\cP(H\Zmp,Y)$.  Next, recall that $[H\Zmp,W]=0$ for each finite
 spectra $W$.  Ravenel proves this in~\cite{ra:lrc} by showing that
 $H\Zmp$ is $E$-acyclic (dissonant) and that finite spectra are
 $E$-local (harmonic), where $E=\bigWedge_{\!p,n}K(n)$.  It was also
 proved earlier by Margolis in~\cite{ma:ems} and by Lin in~\cite{li:dems}
 using the Adams spectral sequence,
 and can be found in Margolis's book~\cite[Cor. 16.27]{ma:ssa}.
 If $\cP(H\Zmp,Y)$ were zero for all $Y$, then the following
 proposition would imply that $H\Zmp=0$, a contradiction.
\end{proof}

\begin{prop}
 If $X$ is $\cA$-projective and $[X,W]=0$ for each $W\in\cF$ then
 $X=0$.  
\end{prop}
\begin{proof}
 We can write $X$ as a retract of a wedge of finite spectra:
 \[ \xy \xymatrix{
     \bigWedge_{\!\alpha} X_\alpha \ar@/l1em/[d] \\
     X \period \ar@/r1em/[u]_i }
 \endxy \]
 Consider $X \xra{i} \bigWedge X_\alpha \xra{j} \prod X_\alpha$,
 where $j$ is the natural map.  
 As $X$ has no maps to finite spectra, the composite $ji$
 is zero.  But
 $\pi_*(j)\:\bigoplus_\alpha\pi_*X_\alpha\ra\prod_\alpha\pi_*X_\alpha$
 is monic, so we see that $\pi_*(i)=0$.  As $i$ is a split
 monomorphism, we know that $\pi_*(X)$ is the image of $\pi_*(i)$.  It
 follows that $X=0$.
\end{proof}

\section{Margolis's axiomatisation conjecture}\label{se:mac}

The Spanier--Whitehead category $\cF$ of finite spectra (in the
classical, topological sense) can be constructed quite simply.
However, all known constructions of the homotopy category $\cS$ of all
spectra are rather intricate.  Moreover, there are a number of
apparently different constructions of this category, all giving the
same result up to equivalence.  (In this section, equivalences of
categories are tacitly required to preserve triangulations and
symmetric monoidal structures.)  It is thus natural to look for a
system of axioms that characterises $\cS$ uniquely in terms of $\cF$.
Margolis~\cite{ma:ssa} conjectured such a characterisation, which
translates into our language as follows: if $\cS'$ is a monogenic
Brown category whose subcategory $\cF'$ of finite objects is
equivalent to $\cF$, then $\cS'$ is equivalent to $\cS$.  As a first
approximation to this conjecture, Margolis showed that $\cS'/\cP'$ is
equivalent to the category $\cH$ of homology theories on $\cF$, or
equivalently to $\cS/\cP$.  Of course, Note~\ref{no:sp} is just a
generalisation of this.  

We can now come somewhat closer to a proof of Margolis's conjecture.
To explain this, we recall some of the theory of linear extensions of
categories.  Our treatment is inspired by~\cite{ba:ah}, but is
different in detail as we only consider additive categories.  
Let $\cB$ be an additive category.  A \dfn{bimodule} over $\cB$ consists
of Abelian groups $D(A,B)$ (for every pair of objects $A,B$ in $\cB$)
together with a trilinear composition operation
\[ \cB(A,B) \otimes D(B,C) \otimes \cB(C,D) \ra D(A,D) \]
written
\[ f\otimes u \otimes g \mapsto f^*g_*u = g_*f^*u. \]
This operation is supposed to have the obvious functoriality
properties.  As an example, because the composite of two phantom maps is
trivial, there is a well-defined composition
\[ \cSP(A,B) \otimes \cP(B,C) \otimes \cSP(C,D) \ra \cP(A,D). \]
This makes $\cP$ into a bimodule over $\cSP$.

If we have an additive functor $F\:\cA\ra\cB$ and a bimodule $D$ over
$\cB$, then we can define a bimodule $F^*D$ over $\cA$ by
$F^*D(A,B)=D(FA,FB)$.  If $F$ is naturally isomorphic to $G$ then one
can check that $F^*D$ and $G^*D$ are isomorphic as bimodules.

A \dfn{linear extension} of $\cB$ by a bimodule $D$ is a category
$\cC$ with the same objects as $\cB$, together with short exact
sequences 
\[ D(A,B) \xra{j} \cC(A,B) \xra{p} \cB(A,B) \]
such that $p$ is a functor and $j(p(f)^*p(g)_*u)=g\circ j(u)\circ f$.
Two such extensions are considered equivalent if there is a functor
$\epsilon\:\cC\ra\cC'$ with $p'\epsilon=p$ and $\epsilon j=j'$ (strict
equalities of functors, not just natural isomorphisms).  We write
$M(\cB,D)$ for the collection of equivalence classes of linear
extensions of $\cB$ by $D$.  The main example of interest to us is of
course the extension $\cP\ra\cS\ra\cSP$.

Suppose again that we have an additive functor $F\:\cA\ra\cB$ and a
linear extension $D\ra\cC\ra\cB$.  Given objects $A,B$ in $\cA$ we
define $F^*\cC(A,B)$ by the pullback diagram
\[ \xy \xymatrix{
  F^*\cC(A,B) \ar[d] \ar[r] & \cC(FA,FB) \ar[d]^p \\
  \cA(A,B)    \ar[r]_F      & \cB(FA,FB) \period
  } \endxy 
\]
One can check that $F^*\cC$ becomes a linear extension of $\cA$ by
$F^*D$.  Moreover, if $G$ is naturally isomorphic to $F$ then $G^*\cC$
is equivalent to $F^*\cC$ as a linear extension.  Thus, a natural
equivalence class of functors $\cA\ra\cB$ induces a map 
$M(\cB,D)\ra M(\cA,F^*D)$.  It is clear that this is essentially
functorial, and thus $M(\cB,D)\simeq M(\cA,F^*D)$ if $F$ is an
equivalence of categories.

A procedure analogous to the Baer sum of extensions makes $M(\cB,D)$
into an Abelian group.  For any pair of objects $A,B$ in $\cB$, the
evident map $M(\cB,D)\ra\Ext(\cB(A,B),D(A,B))$ is a homomorphism.
Unfortunately, this is almost all the information that we have about
the group $M(\cB,D)$ in the cases of interest.  We do not even know
whether $M(\cB,D)$ is a set or a proper class.

We now return to the context of the Margolis conjecture.  We have an
equivalence $F\:\cSP\simeq\cS'/\cP'$.  It follows from
Theorem~\ref{th:pff} that there is a canonical equivalence
$\cP\simeq F^*\cP'$ of bimodules over $\cSP$.  Thus, Margolis's
conjecture is true up to an extension problem.  Together with $F$,
the above equivalence induces a canonical isomorphism
$M(\cS'/\cP',\cP')\simeq M(\cSP,\cP)$.  We need to know whether the
class $u(\cS')$ in $M(\cS'/\cP',\cP')$ that classifies the extension
$\cP'\ra\cS'\ra\cS'/\cP'$ maps to the analogous class 
$u(\cS)\in M(\cSP,\cP)$.  This would follow from Margolis's conjecture.
Conversely, it would almost imply the conjecture, apart from possible
questions about preservation of the triangulation and the monoidal
structure. 

We shall show in the next section that for each $p$ we can choose
spectra $A$ and $B$ such that the image of $u(\cS)$ in
$\Ext(\cSP(A,B),\cP(A,B))$ is not divisible by $p$, and is not
annihilated by any integer $n>0$.
It follows that the same is true of $u(\cS)$ itself.
In particular, we will see that $u(\cS)$ is non-zero.  
This implies that there is no functorial way to choose a representing 
spectrum for a homology theory.

\section{Phantom cohomology}\label{se:pc}

In this section we restrict attention to the classical stable homotopy
category; a more axiomatic approach would yield only a small amount of
extra generality.  Recall that for each Abelian group $A$ there is an
essentially unique spectrum $HA$ with $\pi_0HA=A$ and $\pi_kHA=0$ for
all $k\neq 0$, and that $[X,HA]=H^0(X;A)$.  These objects are called
Eilenberg--Mac\,Lane spectra.  We shall study phantom cohomology classes,
in other words, phantom maps from arbitrary spectra to
Eilenberg--Mac\,Lane spectra.

We start with some algebraic preliminaries.

\begin{defn}
 A monomorphism $B\ra C$ of Abelian groups is said to be \dfn{pure} if
 for each $n > 0$ the induced map $B/n \ra C/n$ is monic.  If we
 regard $B$ as a subgroup of $C$, this says that $nC=(nB)\cap C$.  A
 short exact sequence $B \ra C \ra A$ is said to be \dfn{pure} if the
 map $B \ra C$ is.
\end{defn}

Let $B \ra C \ra A$ be a short exact sequence.
The six term exact sequence involving $\Hom(\Zmn,-)$ and 
$\Ext(\Zmn,-)$ reads
\[ 0 \ra {}_nB \ra {}_nC \ra {}_nA \ra B/n \ra C/n \ra A/n \ra 0 , \]
where we use the notation
\[ {}_nA := \{a \in A \st na = 0 \} \]
and the identifications ${}_nA = \Hom(\Zmn,A)$
and $A/n = \Ext(\Zmn,A)$.
Thus it is clear that pureness of the short exact sequence 
is equivalent to the requirement that ${}_nA \ra B/n$ be zero, or
that ${}_nC \ra {}_nA$ be epic.

Now we present an algebraic proposition which summarises results
that can be found, for example, in~\cite{fu:iag}.

\begin{prop}
 Consider an element $u\in\Ext(A,B)$, corresponding to an extension
 $B\ra C\ra A$.  
 The following are equivalent:
 \begin{itemize}
 \item[(a)] The extension is pure.
 \item[(b)] For each map $A' \ra A$ with $A'$ finitely generated,
   the image of $u$ in $\Ext(A',B)$ is zero.
 \item[(c)] For each $n > 0$, $u \in n \Ext(A,B)$.
   That is, $u$ is in $\bigcap_n n \Ext(A,B)$, the first Ulm subgroup of
   $\Ext(A,B)$.
 \item[(d)] For each map $B \ra B'$ with $B'$ finite,
   the image of $u$ in $\Ext(A,B')$ is zero.
 \end{itemize}
\end{prop}

We define the \dfn{phantom Ext group} $\PExt(A,B)$ to be the subgroup
of $\Ext(A,B)$ consisting of all elements $u$ satisfying the above
conditions.  
These are the phantom maps from $A$ to $\Sigma B$ in $D(\Z)$, 
the derived category of the integers.
It is easy to see that $\PExt$ is a subfunctor of $\Ext$.

\begin{proof}
 (a)$\Rightarrow$(b):
  If suffices to prove (b) when $A' = \Zmn$, as any finitely
  generated group is a sum of cyclic groups, and $\Z$ is projective.
  Given a map $f : \Zmn \ra A$, the class $f^*u$ in $\Ext(\Zmn,B)$
  is zero if and only if $f$ factors through $C \ra A$.
  Now $f$ corresponds to an element of ${}_nA$, and since we
  are assuming $u$ is pure, we know that ${}_nC \ra {}_nA$ is epic
  and can therefore factor $f$ through $C$.  Thus $f^*u = 0$.

 (b)$\Rightarrow$(c):
  We will show that $u$ is in the image of the endomorphism of
  $\Ext(A,B)$ induced by $n : A \ra A$.
  Consider the inclusion $i : {}_nA \ra A$.
  By \cite[Lemma~17.2]{fu:iag}, a bounded torsion group is a sum
  of cyclic groups.  Thus $i^* u = 0$.
  Now in the diagram
  \[ \xy \xymatrix{
    \Ext(A,B) \ar[dr]^n \ar@{->>}[d] \\
    \Ext(nA,B) \ar[r] & \Ext(A,B) \ar[r]^{i^*} & \Ext({}_nA,B) \comma
    } \endxy 
  \]
  the row is exact and the vertical map is an epimorphism
  (because $\Ext^2 = 0$), so $u$ is in the image of multiplication by $n$.

 (c)$\Rightarrow$(d):
  Let $f : B \ra B'$ be a map with $B'$ finite.  To see that the
  image of $u$ in $\Ext(A,B')$ is zero, it suffices to check
  this when $B'$ is $\Zmn$, since a finite group is a product
  of finite cyclic groups.
  But $n$ kills $\Ext(A,\Zmn)$, so if $u$ is a multiple of $n$,
  then $f_*u = 0$.

 (d)$\Rightarrow$(a):
 Finally, assume that for any map $B \ra B'$ with $B'$ finite, the
 image of $u$ in $\Ext(A,B')$ is zero.  Choose an element $b \in B$
 with $b \not\in nB$.  We will show $b \not\in n C$.  (For notational
 simplicity we regard $B$ as a subgroup of $C$.)  Let $K$ be a maximal
 subgroup of $B$ containing $nB$ but not $b$.  The quotient $B/K$ can
 be shown to be ``cocyclic'' and so by \cite[Section~3]{fu:iag} $B/K$
 is isomorphic to $\Zmpk$ for some prime $p$ and some $k$ with
 $1 \leq k \leq \infty$.  Therefore, by assumption (for finite $k$) or
 since $\Zmpi$ is divisible, the quotient map $B \ra B/K$ extends over
 $B \ra C$.  By the choice of $K$, the image of $b$ in $B/K$ is
 non-zero, but the image of $n C$ is zero, so $b \not\in n C$.
\end{proof}

\begin{note}\label{no:pext}
 Clearly, if $A$ is finitely generated, or if there is an integer $n$
 such that $nA=0$, then $\PExt(A,B)=0$ for all $B$.  

 More generally, if $A$ is a torsion group then $A=\dirlim_n\,{}_{n!}A$
 so there is a short exact sequence
 \[ \bigoplus_n {}_{n!}A \ra \bigoplus_n {}_{n!}A \ra A \]
 and a resulting short exact sequence
 \[ \invlim^1\Hom({}_{n!}A,B) \ra \Ext(A,B) \ra 
    \invlim \Ext({}_{n!}A,B).
 \]
 Using part~(b) of the definition of $\PExt(A,B)$, we see that
 \[ \PExt(A,B)=\invlim^1_n \Hom({}_{n!}A,B). \]
\end{note}

Our reason for introducing the phantom Ext groups is the following
theorem, in which $H$ denotes the integral Eilenberg--Mac\,Lane
spectrum.  See the introduction for references to more general
results.

\begin{thm}\label{th:pxhb}
 For any spectrum $X$ and Abelian group $B$ we have
 $\cP(X,HB)=\PExt(H_{-1}X,B)$.
\end{thm}

\begin{proof}
 We begin by describing a map $\cP(X,HB) \ra \PExt(H_{-1}X,B)$.
 Let $u\:X\ra HB$ be a phantom map.  If $Y$ is the cofibre of $u$,
 then we have a short exact sequence $0 \ra B \ra H_0 Y \ra H_{-1} X \ra 0$,
 since $H_0(HB) = B$ by the Hurewicz theorem and since $H_*(u) = 0$
 by Proposition~\ref{pr:ph1}.  We claim that this is a phantom extension,
 and we prove this by showing that for each $n$ the map
 ${}_n(H_0 Y) \ra {}_n(H_{-1}X)$ is surjective.
 Let $a$ be an element of $H_{-1} X$ with $na = 0$.
 This corresponds to a map $S^{-1}\ra H\Smash X$ which can be
 extended to give a map $a'\:S^{-1}/n\ra H\Smash X$.  
 Since phantoms form an ideal under the smash product,
 the composite $S^{-1}/n \ra H \Smash X \ra H \Smash HB$ is
 null and $a'$ factors through $H \Smash Y$.
 Thus $S^{-1} \ra S^{-1}/n \ra H \Smash Y$ represents a class
 in ${}_n(H_0 Y)$ mapping to $a$.

 Conversely, consider the composite
 \[ \PExt(H_{-1}X,B) \ra \Ext(H_{-1}X,B) \ra [X,HB]=H^0(X;B) , \]
 where the first map is the inclusion and the second map comes from
 the universal coefficient sequence.
 We claim that a map $u$ in the image of this composite is a phantom map.
 Indeed, if $W$ is a finite spectrum and $W \ra X$ is a map, 
 then by naturality the restriction of $u$ to $W$ lies in the image of
 $\PExt(H_{-1}W,B)$, which is trivial because $H_{-1}W$ is finitely
 generated.  It follows that $u$ is a phantom map.

 We leave it to the reader to check that the two maps we have
 constructed are inverses.
\end{proof}

This allows us to calculate all phantom maps between
Eilenberg--Mac\,Lane spectra.

\begin{cor}
 We have
 \[ \cP(\Sigma^k HA,HB) = \begin{cases}
     \PExt(A,B) & \text{ if } k = -1 \\
     0          & \text{ otherwise. }
    \end{cases}
 \]
\end{cor}
\begin{proof}
 For $j<0$ we have $H_jHA=0$, and $H_0HA=A$ by the Hurewicz theorem.
 Given this, the claim follows for $k\geq -1$ by a simple application
 of Theorem~\ref{th:pxhb}.  For $j>0$ we may have $H_jHA\neq 0$, but
 we claim that $\PExt(H_jHA,B)=0$ nonetheless; this will cover the
 case $k<-1$.  To see this, fix $j>0$ and let $\{A_\alpha\}$ be the
 directed set of finitely generated subgroups of $A$.  
 The natural map $\dirlim_\alpha (HA_\alpha)_* X \ra (HA)_* X$
 is an isomorphism for each $X$, since it is when $X$ is a sphere,
 and both sides are homology theories.
 Taking $X = H$ we find that
 $H_jHA=\dirlim_\alpha H_jHA_\alpha$.  By working rationally, we see
 that $H_jHA$ is a torsion group, so it is the direct sum of its
 localisations at different primes.  We claim that
 $H_j(HA_\alpha)_{(p)}$ is a vector space over $\Zmp$.  Using the fact
 that $H_jHA_\alpha=(HA_\alpha)_jH$ and the fact that the universal
 coefficient sequence splits, we are reduced to proving that $H_jH$ is
 killed by $p$.  This is a classical calculation; an account appears
 in~\cite{ko:ico}.  This implies that $H_jHA$ is a direct sum of
 (prime) cyclic groups; it follows easily that $\PExt(H_jHA,B)=0$ as
 required. 
\end{proof}

We next study a special case in which the short exact sequence
\[ \cP(\Sigma^{-1}HA,HB) \ra \cS(\Sigma^{-1}HA,HB) \ra 
        \cA(\Sigma^{-1}HA,HB)
\]
can be understood explicitly.  We choose a prime $p$ and take
\[ A = \Zmpi = \Q/\Z_{(p)} = \dirlim_k \Zmpk. \]
For the moment we consider an arbitrary Abelian group $B$.  As in
Note~\ref{no:pext}, we have a short exact sequence
\[ \PExt(A,B) \ra \Ext(A,B) \ra \invlim_k\Ext(\Zmpk,B) . \] 
Note that $\Ext(\Zmpk,B)=B/p^k$, so the third term is just the
$p$-completion $\Bh$ of $B$.  The middle term is the
Ext-$p$-completion of $B$, as studied in~\cite{boka:hlc}; we shall
denote it by $\Bt$.  And it is clear that the first term is
\[ \pinf \Bt = \bigcap_k p^k\Bt , \]
since everything is $p$-local.
Using the fact that $\cS(\Sigma^{-1} HA,HB) = \Ext(A,B)$, we
find that our phantom sequence is just
\[ \pinf \Bt \ra \Bt \ra \Bh . \]

It is tempting to believe that $\pinf\Bt$ is a divisible group, but
this is never true unless $\pinf\Bt=0$.  Any element of $\pinf\Bt$ is
divisible by $p$ in $\Bt$ but not necessarily in $\pinf\Bt$.

Let $B$ be a free Abelian group, say $B=\bigoplus_{k=0}^\infty\Z$.
Then $\Hom(\Zmpj,B)=0$ so $\pinf\Bt=\invlim^1_j\Hom(\Zmpj,B)=0$ so
$\Bt=\Bh$.  Let $v(a)$ denote the $p$-adic valuation of a $p$-adic
integer $a\in\Zp$.  It is not hard to see that
\[ \Ext(\Zmpi,B)=\Bt=\Bh=\{\ab\in\prod_k\Zp\st v(a_k)\ra\infty\}. \]
One can also see directly that $\Hom(\Zmpi,B)=0$.

Now consider the case $B=\bigoplus_k\Zmpk$.  We then have a short
exact sequence
\[ \bigoplus_k\Z \xra{f} \bigoplus_k\Z \ra B \]
where $f$ is multiplication by $p^k$ on the $k$'th factor.  One can
again see directly that $\Hom(\Zmpi,B)=0$.  The six-term exact
sequence obtained by applying the functors $\Hom(\Zmpi,-)$ and
$\Ext(\Zmpi,-)$ to the above presentation of $B$ therefore collapses
to a short exact sequence
\[ \Ext(\Zmpi,\bigoplus_k\Z)\xra{f}\Ext(\Zmpi,\bigoplus_k\Z)\ra\Bt.\]
It follows using the previous paragraph that
\[ \Bt = 
   \{\ab\st v(a_k)\ra\infty\}/\{\ab\st 0\le v(a_k)-k\ra\infty\}.
\]
One can also see directly that 
\[ \Bh = \{\ab\st v(a_k)\ra\infty\}/\{\ab\st 0\le v(a_k)-k\}. \]
It follows that $\pinf\Bt$ (which is the kernel of the map
$\Bt\ra\Bh$) is given by
\[ \pinf\Bt =
    \{\ab\st 0\le v(a_k)-k\}/\{\ab\st 0\le v(a_k)-k\ra\infty\},
\]
and this can also be expressed as 
\[ \pinf\Bt =
     \prod_k\Zp/\biggl\{\bb\in\prod_k\Zp\st v(b_k)\ra\infty\biggr\}
\]
(where $a_k=p^kb_k$).  It is easy to see from this that $\pinf\Bt$ is
nonzero and torsion-free.

We now return to the case of a general Abelian group $B$.  Let
$w\in\Ext(\Bh,\pinf\Bt)$ be the element classifying the canonical
sequence $\pinf\Bt\ra\Bt\ra\Bh$, and let
\[ \delta\:\Hom(\Zmp,\Bh)\ra\Ext(\Zmp,\pinf\Bt) \]
be the obvious connecting homomorphism.
 
\begin{prop}\label{pr:wbi}
 Let $B$ be an Abelian group.  The following are equivalent:
 \begin{roenumerate}
  \item $\pinf\Bt=0$.
  \item The natural map $\Bt\ra\Bh$ is an isomorphism.
  \item $w=0$.
  \item $w$ is divisible by $p$.
  \item $\delta=0$.
  \item $\delta$ is divisible by $p$.
 \end{roenumerate}
\end{prop}

\begin{proof}
 (i)$\Rightarrow$(ii)$\Rightarrow$(iii)
 $\Rightarrow$(iv)$\Rightarrow$(vi): easy.

 (vi)$\Rightarrow$(v): This is also clear, as the source and target of
 $\delta$ are killed by $p$.

 (v)$\Rightarrow$(i): The next map in the sequence is
 $\Ext(\Zmp,\pinf\Bt)\ra\Ext(\Zmp,\Bt)$, which can be identified with
 the natural map $(\pinf\Bt)\quo p\ra\Bt\quo p$.  But this latter map
 is clearly zero, so the connecting homomorphism $\delta$ is epic.
 Its image is $(\pinf\Bt)\quo p$, so this group is zero, so $\pinf\Bt$
 is $p$-divisible.  This means that $\pinf\Bt=\Hom(\Z,\pinf\Bt)$ is a
 quotient of $\Hom(\Z[\pinv],\pinf\Bt)$, which is a subgroup of
 $\Hom(\Z[\pinv],\Bt)$.  However, \cite[VI.3.4]{boka:hlc} tells us
 that $\Hom(\Z[\pinv],\Bt)=0$; it follows that $\pinf\Bt=0$.
\end{proof}

If $B=\bigoplus_k\Zmpk$, then we have $\pinf\Bt\neq 0$ and thus $w$
is not divisible by $p$.  Our next result will show that $w$ has
infinite order.

\begin{prop}\label{pr:wbii}
 Let $B$ be an Abelian group such that $p^kw=0$.  Then $p^k\pinf\Bt=0$.  
\end{prop}
\begin{proof}
 Let $i\:\pinf\Bt\ra\Bt$ and $q\:\Bt\ra\Bh$ be the usual maps.  Let
 $C$ be the pullback of $\Bt$ along the map $p^k\:\Bh\ra\Bh$, so
 $C=\{(a,b)\in\Bt\times\Bh\st q(a)=p^kb\}$.  The hypothesis
 $p^kw=0$ means that the evident sequence $\pinf\Bt\ra C\ra\Bh$ is
 split; the splitting map $\Bh\ra C$ necessarily has the form
 $c\mapsto(f(c),c)$, where $qf=p^k\:\Bh\ra\Bh$.  The functor
 $A\mapsto\pinf A$ preserves split exact sequences, and $\pinf\Bh=0$
 (directly from the definitions) so $\pinf C=\pinf\pinf\Bt$.  On the
 other hand, suppose that $b\in\pinf\Bt$, say $b=p^ib_i$ for each $i$,
 with $b_i\in\Bt$.  Then $(p^kb_i,q(b_i))\in C$ and
 $p^i(p^kb_i,q(b_i))=(p^kb,0)$.  It follows that
 $p^k\pinf\Bt\leq\pinf C=\pinf\pinf\Bt$.  This means that
 $p^k\pinf\Bt$ is a divisible subgroup of $\Bt$; as in the proof of
 Proposition~\ref{pr:wbi}, we conclude that $p^k\pinf\Bt=0$.
\end{proof}

Now take $B=\bigoplus_k\Zmpk$ again.  We saw previously that
$\pinf\Bt$ is non-trivial and torsion-free.  It follows easily that
$p^kw\neq 0$ for all $k$.

We can now prove a result stated in Section~\ref{se:mac}.  Recall that
we defined there a group $M(\cSP,\cP)$ and an element 
$u\in M(\cSP,\cP)$ that classifies the linear extension of categories
$\cP\ra\cS\ra\cSP$.  The image of $u$ under a certain homomorphism
$M(\cSP,\cP)\ra\Ext(\Bh,\pinf\Bt)$ is $w$.  It follows that $u$ is
not divisible by $p$ for any prime, so $u$ is not divisible by any
integer $n>1$.  If $u$ were annihilated by any $m>0$ then the image of
$u$ in any $p$-local group (such as $\Ext(\Bh,\pinf\Bt)$) would be
annihilated by some power of $p$.  Thus, we conclude that $u$ does not
have finite order.

\section{Universal homology theories}\label{se:uht}

Although one is mostly interested in homology theories with values in
the category of Abelian groups, one can also consider more general
Abelian categories.  In this section, we recall a construction of
Freyd~\cite{fr:sh} which gives a universal example of an Abelian
category $\cB$ equipped with a homology theory $\cS\ra\cB$.  We also
show that the functor $\yoneda\:\cS\ra\cA$ is the universal example of
a homology theory with values in an Abelian category satisfying
Grothendieck's axiom AB 5.

At one point we need a fact that holds in all monogenic Brown
categories that we care about, but which we have not been able to
deduce from the axioms (although we suspect that it does follow).  For
simplicity, we therefore restrict attention to the classical stable
homotopy category.

Let $\cB$ be the following category.  The objects of $\cB$ are just
the morphisms of $\cS$.  Given a map $u\:W\ra X$ in $\cS$, we shall
write $I(u)$ for $u$ thought of as an object of $\cB$.  The group
$\cB(I(u),I(v))$ is the quotient of the group of commutative squares
\[ \xy \xymatrix{
  W \ar[d]_u \ar[r]^f & Y \ar[d]^v \\
  X \ar[r]_g          & Z
  } \endxy 
\]
by the subgroup of squares for which the map $vf=gu$ vanishes.  This
gives a category $\cB$ in an obvious way.  There is a full and
faithful embedding $J\:\cS\ra\cB$ sending $X$ to $I(1_X)$.  Freyd
shows that $\cB$ is an Abelian category and that $J$ is a homology
theory.  Given a morphism $u\:W\ra X$ in $\cS$, the image of the
morphism $Ju\:J(W)\ra J(X)$ is just $I(u)$.  Moreover, the image of
$J$ is the subcategory of injective objects in $\cB$, which is the
same as the subcategory of projective objects.  Freyd also shows that
for any Abelian category $\cC$ and any homology theory $K\:\cS\ra\cC$,
there is an essentially unique strongly additive exact functor
$K'\:\cB\ra\cC$ such that $K'J\simeq K$.  (We say that a functor is
\dfn{strongly additive} if it preserves all coproducts.  Freyd
actually proves the corresponding result without strong additivity
but the necessary modifications are trivial.)  In fact, $K'I(u)$ is
just the image of the morphism $Ku$ in $\cC$.  In particular, this
construction gives a functor $\cB\ra\cA$.

The following result is analogous to Theorem~\ref{th:hpi}.
\begin{prop}
 In $\cB$, any object of finite projective or injective dimension is
 both projective and injective, and thus lies in the image of $J$.
\end{prop}
\begin{proof}
 Suppose that $X$ has projective dimension at most $n>0$.  There is
 then a short exact sequence $Y\ra P\ra X$, where $Y$ has projective
 dimension at most $n-1$; by induction, we may assume that $Y$ is
 projective.  As projectives are injective, the sequence splits, so
 $X$ is a retract of $P$ and thus is projective.
\end{proof}

While this seems a pleasant construction, the finiteness properties of
the category $\cB$ are poor.  We believe that every nonzero object
has a proper class of subobjects, for example.

Next, recall that an Abelian category is said to satisfy AB 5 if
set-indexed colimits exist and filtered colimits are
exact~\cite{gr:qpa}.  The category of Abelian groups satisfies
AB 5, as does the functor category $\cA$.  

\begin{prop}\label{pr:abfiveh}
 If $\cC$ is an Abelian category satisfying AB 5 and $K\:\cS\ra\cC$ is
 a homology theory, then $KX=\dirlim_{\Lambda(X)}KX_\alpha$.  Thus,
 $Kf$ is zero for any phantom map $f$.
\end{prop}
\begin{proof}
 Define $\widehat{K}X=\dirlim_{\Lambda(X)}KX_\alpha$.  Here we will
 need to use the fact (mentioned after Corollary~\ref{cor:lambda})
 that $\Lambda(X)$ is the diagram of all pairs $(U,u)$, where $U$ lies
 in some small skeleton of $\cF$ and $u\:U\ra X$.  Using this we see
 that $\widehat{K}$ is an additive functor $\cS\ra\cC$, and that there
 is an evident natural map $\widehat{K}\ra K$
 (compare~\cite[Proposition 2.3.9]{hopast:ash}).  If $X$ is finite
 then $\Lambda(X)$ has a terminal object, so that $\widehat{K}X=KX$.
 If we can show that $\widehat{K}$ preserves coproducts and sends
 cofibre sequences to exact sequences, then the usual argument will show
 that $\widehat{K}X=KX$ for all $X$.  Consider a cofibre sequence
 $X\ra Y\ra Z$.  We may assume that $X$ is a CW subspectrum of $Y$,
 and that $Z$ is the quotient.  Let $\{Y_\alpha\st\alpha\in I\}$ be
 the directed set of finite subspectra of $Y$.  Write
 $X_\alpha=Y_\alpha\cap X$ and $Z_\alpha=Y_\alpha/X_\alpha$, so we
 have a cofibre sequence $X_\alpha\ra Y_\alpha\ra Z_\alpha$ for each
 $\alpha$.  It is easy to see that the evident functors from $I$ to
 $\Lambda(X)$, $\Lambda(Y)$ and $\Lambda(Z)$ are cofinal, so that
 $\widehat{K}X=\dirlim_IKX_\alpha$ and so on.  As direct limits are
 exact, we conclude that the sequence
 $\widehat{K}X\ra\widehat{K}Y\ra\widehat{K}Z$ is exact as required. 
 
 We next verify that $\widehat{K}$ preserves coproducts.  Consider a
 family of spectra $\{X_i\st i\in I\}$.  Let $\Lambda$ be the full
 subcategory of $\prod_I\Lambda(X_i)$ consisting of those objects
 $(Z_i)_{i\in I}$ such that $Z_i=0$ for almost all $i$.  It is not
 hard to see that this is a filtered category, and that the
 projections $\Lambda\ra\Lambda(X_i)$ are cofinal functors.  The
 functor from $\Lambda$ to $\Lambda(\bigWedge_{\!i} X_i)$ is also
 cofinal.  By writing $\widehat{K}(X_i)$ and
 $\widehat{K}(\bigWedge_{\!i} X_i)$ as colimits indexed by $\Lambda$,
 we see that
 $\widehat{K}(\bigWedge_{\!i}X_i)=\bigoplus_i \widehat{K}(X_i)$ as
 required.
\end{proof}

In a more general monogenic Brown category, it is more difficult to
prove that $\widehat{K}$ is an exact functor.  An obvious approach is
to replace the sequences $X_\alpha \ra Y_\alpha \ra Z_\alpha$
considered above by the category of
all cofibre sequences of finite objects equipped with a map to the
sequence $X\ra Y\ra Z$.  However, it is not clear that this is a
filtered category.  The difficulty is related to the existence of maps
of cofibre sequences that are not good in the sense of
Neeman~\cite{ne:nat}.

We can deduce from the above that $\yoneda\:\cS\ra\cA$ is the
universal example of a homology theory with values in an AB 5
category.

\begin{prop}
 Let $\cC$ be an AB 5 category, and $K\:\cS\ra\cC$ a homology theory.
 Then there is an essentially unique strongly additive exact functor
 $K'\:\cA\ra\cC$ such that $K'\circ\yoneda\simeq K$.
\end{prop}

\begin{proof}
 First, let $\cH$ be the category of homology theories, so
 $\cH\subset\cA$ and $\cH\simeq\cS/\cP$.  By
 Proposition~\ref{pr:abfiveh}, we know that $K$ kills phantom maps, so
 it factors in an essentially unique way through
 $\yoneda\:\cS\ra\cH$.  We write $K$ again for the resulting functor
 $\cH\ra\cC$.  As the cofibre of an $\cA$-epimorphism is phantom, we
 see that $K$ sends epimorphisms of homology theories to epimorphisms,
 and similarly for monomorphisms.
  
 Consider an object $F\in\cA$.  We know that $\cA$ has enough
 projectives and injectives, so we can choose maps $P\xra{f}F\xra{g}I$
 where $f$ is epic, $g$ is monic, $P$ is projective and $I$ is
 injective.  In particular, $P$ and $I$ are homology theories, so
 $K(P)$ and $K(I)$ are defined.  We would like to define $K'(F)$ to be
 the image of the map $K(gf)\:K(P)\ra K(I)$; we need only check that
 this is well-defined.  Indeed, if we chose a different epimorphism
 $f'\:P'\ra F$ then we could use the projectivity of $P$ and $P'$ to
 show that $f$ and $f'$ factor through each other; it follows easily
 that $K(gf)$ and $K(gf')$ have the same image, regarded as a
 subobject of $K(I)$.  A similar argument shows that our definition is
 essentially independent of $g$.

 Next, consider a morphism $v\:F\ra G$ in $\cA$.  Choose sequences
 $P\ra F\ra I$ and $Q\ra G\ra J$ as above.  Using the projectivity
 of $P$ and the injectivity of $J$, we can choose maps $u\:P\ra Q$
 and $w\:I\ra J$ compatible with $v$.  These induce a map 
 $K(F)\ra K(G)$, which we would like to call $K'(v)$.  We must check
 that this does not depend on the choice of $u$ and $w$.  An easy
 argument reduces us to the case $v=0$; this implies that the diagonal
 map in the square 
 \[ \xy \xymatrix{
  P \ar[d] \ar[r]^u & Q \ar[d] \\
  I \ar[r]_w        & J
  } \endxy 
 \]
 is zero, and thus the induced map
 $\image(K(P\ra I))\ra\image(K(Q\ra J))$ is zero as required.  Our
 definition of $K'(v)$ is thus unambiguous, and it is easy to see that
 it gives a functor. 

 Suppose that $F$ is a homology theory.  Then $P\ra F$ is an
 epimorphism of homology theories, so $K(P)\ra K(F)$ is epic.
 Similarly, $K(F)\ra K(I)$ is monic.  It follows directly that
 $K'(F)=K(F)$.  Thus, $K'$ is an extension of $K$.

 If $v\:F\ra G$ is a monomorphism then we may choose $I=J$ and $w=1$;
 this makes it clear that $K'(v)$ is a monomorphism.  Similarly, $K'$
 preserves epimorphisms.

 We next show that $K'$ preserves kernels.  Consider a map
 $v\:F\ra G$.  Choose an epimorphism $f\:P\ra F$ and a monomorphism
 $g\:G\ra J$.  As $P$ and $J$ are homology theories, we can choose a
 map of spectra inducing the map $gvf\:P\ra J$, and let $j\:H\ra P$ be
 its fibre.  As $H\ra P\ra J$ is zero, we see that $H\ra P\ra F$
 factors through $\ker(gv)=\ker(v)$.  As $K'$ preserves monomorphisms
 and epimorphisms, we obtain a diagram as follows:
 \[ \xy \xymatrix{
   K'(H) \ar[d] \ar[r]^{K'j}& K'(P) \ar@{->>}[d]_{K'f} \\
   K'(\ker(v))\; \ar@{>->}[r] & K'(F) \ar[r]_{K'v} & K'(G)\; 
    \ar@{>->}[r]_{K'g} & K'(J) \period
   } \endxy
 \]
 As $H\ra P\ra J$ comes from a cofibre sequence of spectra, we know
 that $K'(H)\ra K'(P)\ra K'(J)$ is exact.  A diagram chase (using
 elements in the sense of~\cite{ma:cwm}, for example) now shows that
 $K'(\ker(v))\ra K'(F)\ra K'(G)$ is exact as required.

 Similarly, we see that $K'$ preserves cokernels; it is thus an exact
 functor.  

 Finally, we need to show that $K'$ preserves coproducts.  Consider a
 family $\{F_i\}$ of objects of $\cA$, and choose maps
 $P_i\ra F_i\ra I_i$ in the usual way.  Write $P=\bigoplus_iP_i$ and
 $F=\bigoplus_iF_i$ and $I=\bigoplus_iI_i$, so we have an epimorphism
 $P\ra F$ and a monomorphism $F\ra I$ (but $I$ need not be
 injective).  As $K'$ is exact, we see that $K'(F)$ is the image of
 $K'(P)\ra K'(I)$.  As $K$ preserves coproducts of spectra, we see
 that $K'$ preserves coproducts of homology theories, so
 $K'(P)=\bigoplus_i K'(P_i)$.  Similarly,
 $K'(I)=\bigoplus_i K'(I_i)$.  It follows that
 $K'(F)=\bigoplus_i K'(F_i)$ as required.

 It is also clear that any extension of $K$ that preserves images (in
 particular, any exact extension of $K$) must be equivalent to $K'$.
\end{proof}

We conclude this section with an interesting, if somewhat disconnected
result.  Consider an essentially small additive category $\cF$.
Recall that there is an essentially unique category $\Ind(\cF)$ (the
Ind completion of $\cF$) equipped with a full and faithful embedding
$\cF\ra\Ind(\cF)$ (thought of as an inclusion) such that
\begin{roenumerate}
 \item $\Ind(\cF)$ has colimits for all small filtered diagrams.
 \item Every object of $\Ind(\cF)$ is the colimit of a small filtered
  diagram of objects of $\cF$.
 \item If $X$ is an object of $\cF$ then the functor $\Ind(\cF)(X,-)$
  preserves filtered colimits.
\end{roenumerate}
The Ind completion of a category was introduced in~\cite{gr:td} and
was described in detail in~\cite{grve:sga4}.

This category can be constructed in (at least) two ways.  The first
way is to consider pairs $(I,X)$ where $I$ is a small filtered
category and $X$ is a functor $I\ra\cF$.  We define $\Ind(\cF)$ to be
the category of such pairs, with morphisms
\[ \Ind(\cF)((I,X),(J,Y)) = \invlim_I\dirlim_J\cF(X_i,Y_j). \]
Alternatively, we can embed $\cF$ in the category $\cB$ of additive
functors $\cF^{\textup{op}}\ra\Ab$ by $X\mapsto [-,X]$.  We then
define $\Ind(\cF)$ to be the subcategory of all functors $F\in\cB$
that can be written as a filtered colimit of a small diagram of
objects of $\cF$.  It is equivalent to require that the category of
pairs $(X,a)$ (where $X\in\cF$ and $a\in FX$) is filtered.
 
\begin{thm}
 Let $\cF$ be the category of finite spectra.  Then there is an
 equivalence of categories $\Ind(\cF)=\cH$ \ulp{}where $\cH$ is the
 category of homology theories on $\cF$\urp{}.
\end{thm}
\begin{proof}
 We use the second description of $\Ind(\cF)$, as a subcategory of
 $\cB=[\cF^{\textup{op}},\Ab]$.  Composition with the
 Spanier-Whitehead duality functor gives an equivalence of $\cB$ with
 $\cA=[\cF,\Ab]$, which sends $[-,X]$ to $\yoneda_X$.  Thus,
 $\Ind(\cF)$ is equivalent to the category of those functors
 $\cF\ra\Ab$ that can be written as small filtered colimits of
 functors of the form $\yoneda_X$ where $X$ is small.  As filtered
 colimits are exact, every such functor is a homology theory.
 Conversely, every homology theory is of the form $\yoneda_Y$ for some
 $Y$.  Since $\yoneda_Y=\dirlim_{\Lambda(Y)}\yoneda_{Y_\alpha}$, it
 follows that every homology theory lies in $\Ind(\cF)$.
\end{proof}

We conclude that the Ind completion of a triangulated category need
not be a triangulated category.  For example, consider a monogenic
Brown category with a non-zero phantom map $f\:X \ra Y$.  If $Z$ is
the cofibre of $f$, then the map $\yoneda_Y \ra \yoneda_Z$ is monic
but not split.  However in a triangulated category all monics split,
so $\cH$ is not triangulated.  (Hartshorne mentions in~\cite{ha:rd}
that the Ind completion may not be triangulated, but he does not
indicate a proof.)  We also gain some insight into the Pro completion
of the category of spectra, which has been used by various people for
various purposes, mainly concentrating on towers of spectra rather
than more general inverse systems.  The Pro completion of a category
$\cC$ is just $\Ind(\cC^{\textup{op}})^{\textup{op}}$.  As
$\cF\simeq\cF^{\textup{op}}$ (by Spanier-Whitehead duality) we see
that the Pro category of finite spectra is equivalent to the opposite
of the category of homology theories.  The subcategory consisting of
towers of finite spectra is equivalent to the opposite of the category
of homology theories with countable coefficients.


\end{document}